\documentclass[12pt,a4paper]{article}

\usepackage[left=1.5cm,right=3cm,top=2cm,bottom=2cm]{geometry}
\usepackage{array}
\usepackage[fleqn]{amsmath}
\usepackage{amssymb,latexsym}
\usepackage{amsthm}
\usepackage[section]{placeins}
\usepackage[colorlinks=true,linkcolor=black,citecolor=black,urlcolor=black]{hyperref}
\usepackage{enumitem}
\usepackage{graphicx}
\usepackage{tikz}
\usepackage{todonotes}

\newtheorem{theorem}{Theorem}[section]

\newtheorem{observation}[theorem]{Observation}

\theoremstyle{definition}
\newtheorem{definition}{Definition}

\DeclareMathOperator{\Cay}{Cay}
\DeclareMathOperator{\PSL}{PSL}
\DeclareMathOperator{\SL}{SL}
\DeclareMathOperator{\A}{A}
\newcommand{\Z}{\mathbb{Z}}

\title{Small graphs and hypergraphs of given degree and girth}

\author{
Grahame Erskine\thanks{Open University, Milton Keynes, UK}\footnotemark[1]\\ \texttt{\small grahame.erskine@open.ac.uk}
\and James Tuite\footnotemark[1]\\ \texttt{\small james.tuite@open.ac.uk}
}

\date{}

\begin{document}
\maketitle
\let\thefootnote\relax
\let\thefootnote\relax

\vspace*{-8ex}
\begin{abstract}
\noindent The search for the smallest possible $d$-regular graph of girth $g$ has a long history, and is usually known as the \emph{cage problem}. This problem has a natural extension to hypergraphs, where we may ask for the smallest number of vertices in a $d$-regular, $r$-uniform hypergraph of given (Berge) girth $g$. We show that these two problems are in fact very closely linked. By extending the ideas of Cayley graphs to the hypergraph context, we find smallest known hypergraphs for various parameter sets. Because of the close link to the cage problem from graph theory, we are able to use these techniques to find new record smallest cubic graphs of girths 23, 24, 28, 29, 30, 31 and 32.
\end{abstract}

\section{Introduction}\label{sec:intro}
The problem of finding the smallest possible graph of given girth and minimum degree (the \emph{cage} problem) is an active area of research (see the survey~\cite{exoo}). There is a natural extension of this problem to hypergraphs, and we will be concerned in this paper with both problems and with the links between them. We begin with some elementary definitions.

A \emph{hypergraph} $H$ is a set of vertices $V(H)$ and a set $E(H)$ of subsets of $V(H)$ called \emph{hyperedges} (or simply \emph{edges} if there is no risk of confusion). The \emph{degree} or \emph{valency} of a vertex is the number of hyperedges containing it; if every vertex has the same degree $d$, then we say the hypergraph is $d$-\emph{regular}. If all hyperedges have the same cardinality $r$, then we say the hypergraph is $r$-\emph{uniform}. A \emph{Berge cycle} of length $k$ in a hypergraph is a sequence $v_0,e_0,v_1,e_1,\ldots,v_{k-1},e_{k-1},v_0$ such that each $v_i$ is contained in $e_{i-1}$ and $e_i$ (mod $k$), all $v_i$ are distinct and all $e_i$ are distinct. (Note that other definitions of `cycle' in hypergraphs are possible; in this work we are concerned only with Berge cycles.) The \emph{girth} of a hypergraph is the length of its smallest Berge cycle. We say a hypergraph is \emph{linear} if two distinct hyperedges meet in at most one vertex; thus a hypergraph is linear if and only if its girth is at least 3. With these definitions, a graph may be viewed as a 2-uniform hypergraph, and then the usual graph definitions of degree and girth are consistent with the hypergraph equivalents.

A hypergraph may equivalently be viewed as an incidence structure, and it is natural to associate with the hypergraph $H$ a bipartite incidence graph (or \emph{Levi graph}) $I(H)$ which has vertex set $V(H)$ (black vertices) and $E(H)$ (white vertices); there is an edge from $v$ to $e$ in $I(H)$ if and only if $v\in e$ in $H$. If $H$ is a $d$-regular, $r$-uniform hypergraph, then $I(H)$ is a $(d,r)$-biregular bicoloured graph; we use the term `bicoloured' here to emphasise that the black and white vertices are distinguished. By swapping these colour classes in $I(H)$ we get the incidence graph $I^*$ of the \emph{dual} hypergraph $H^*$. This is $r$-regular and $d$-uniform.

The remainder of the paper is organised as follows. In Section~\ref{sec:duality} we investigate the girth parameter in dual hypergraphs and deduce new lower bounds on the minimum number of vertices of a $d$-regular, $r$-uniform hypergraph in the case where $r>d$. In Section~\ref{sec:cayley} we describe a method of finding small examples of $d$-regular, $r$-uniform hypergraphs of given girth, based on similar methods used in the graph cage problem. In addition, we show how hypergraphs can be used to obtain new examples of small regular graphs of given girth, by exploiting the idea of duality. Finally in Section~\ref{sec:cages} we summarise the results of a large computer search, including a number of new entries in the table of smallest known cubic graphs of given girth.
\section{Dual hypergraphs}\label{sec:duality}
A simple counting argument (see for example Ellis and Linial~\cite{EllLin} following Hoory~\cite{Hoory}) yields the following lower bound $M(d,r,g)$ (the \emph{Moore bound}) on the number of vertices of a $d$-regular, $r$-uniform hypergraph $H$ of girth $g$. If $g=2k+1$ is odd, then
\begin{equation}\label{eq:mooreodd1}
\lvert V(H)\rvert\geq M(d,r,g)=1+d(r-1)\frac{(d-1)^k(r-1)^k-1}{(d-1)(r-1)-1}
\end{equation}
and if $g=2k$ is even, then
\begin{equation}\label{eq:mooreeven1}
\lvert V(H)\rvert\geq M(d,r,g)=r\,\frac{(d-1)^k(r-1)^k-1}{(d-1)(r-1)-1}.
\end{equation}

The following simple observations are fundamental to the remainder of the paper.
\begin{observation}\label{obs:2g}
A (Berge) cycle of length $\ell$ in a hypergraph corresponds to a cycle of length $2\ell$ in its Levi graph.
\end{observation}
\begin{observation}\label{obs:dual}
Since the Levi graph of a hypergraph is isomorphic to that of its dual, the girth of a hypergraph is equal to the girth of its dual.
\end{observation}
Observation~\ref{obs:dual} allows us to obtain a better bound for the minimum order of a $d$-regular, $r$ uniform hypergraph $H$ of girth $2k+1$ than that noted in~\cite{EllLin}, in the case where $r>d$. Since $\lvert V(H)\rvert=\lvert E(H^*)\rvert=\dfrac{r}{d}\lvert V(H^*)\rvert$ and since $H^*$ also has girth $g$, it follows that
\begin{equation}\label{eq:mooreodd2}
\lvert V(H)\rvert\geq \frac{r}{d}\left(1+r(d-1)\frac{(d-1)^k(r-1)^k-1}{(d-1)(r-1)-1}\right).
\end{equation}
Since expression (\ref{eq:mooreodd2}) exceeds expression (\ref{eq:mooreodd1}) if $r>d$, we can find a tighter bound for the minimum order in this case. Since a similar argument holds for hypergraphs of even girth, we have the following result.
\begin{theorem}\label{thm:moorebound}
Let $H$ be a $d$-regular, $r$ uniform hypergraph of girth $g$. If $g=2k+1$ is odd, then
\[
\lvert V(H)\rvert\geq
\begin{cases}
1+d(r-1)\dfrac{(d-1)^k(r-1)^k-1}{(d-1)(r-1)-1} & \text{if } r\leq d;\\[3ex]
\dfrac{r}{d}\left(1+r(d-1)\dfrac{(d-1)^k(r-1)^k-1}{(d-1)(r-1)-1}\right) & \text{if } r>d.
\end{cases}
\]
If $g=2k$ is even, then
\[
\lvert V(H)\rvert\geq
\begin{cases}
r\,\dfrac{(d-1)^k(r-1)^k-1}{(d-1)(r-1)-1} & \text{if } r\leq d;\\[3ex]
\dfrac{r^2}{d}\,\dfrac{(d-1)^k(r-1)^k-1}{(d-1)(r-1)-1} & \text{if } r>d.
\end{cases}
\]
\end{theorem}
The impact of Theorem~\ref{thm:moorebound} can be seen in Table~\ref{tab:newmoore}. 
\begin{table}[h]\centering\small
	\begin{tabular}{|c|ccccccc|ccccccc|}
		\hline
		\multicolumn{1}{|c|}{} & \multicolumn{7}{c|}{Ignoring Theorem~\ref{thm:moorebound}} & \multicolumn{7}{|c|}{Allowing for Theorem~\ref{thm:moorebound}} \\
		$d\setminus r$ & 2 & 3 & 4 & 5 & 6 & 7 & 8 & 2 & 3 & 4 & 5 & 6 & 7 & 8 \\
		\hline
		2 & 5 & 13 & 25 & 41 & 61 & 85 & 113 & 5 & 15 & 34 & 65 & 111 & 175 & 260 \\
		3 & 10 & 31 & 64 & 109 & 166 & 235 & 316 & 10 & 31 & 76 & 152 & 266 & 427 & 643 \\
		4 & 17 & 57 & 121 & 209 & 321 & 457 & 617 & 17 & 57 & 121 & 245 & 434 & 700 & 1058 \\
		5 & 26 & 91 & 196 & 341 & 526 & 751 & 1016 & 26 & 91 & 196 & 341 & 606 & 982 & 1487 \\
		6 & 37 & 133 & 289 & 505 & 781 & 1117 & 1513 & 37 & 133 & 289 & 505 & 781 & 1267 & 1922 \\
		7 & 50 & 183 & 400 & 701 & 1086 & 1555 & 2108 & 50 & 183 & 400 & 701 & 1086 & 1555 & 2360 \\
		8 & 65 & 241 & 529 & 929 & 1441 & 2065 & 2801 & 65 & 241 & 529 & 929 & 1441 & 2065 & 2801 \\
		\hline
	\end{tabular}
	\caption{Order bounds for $d$-regular, $r$-uniform hypergraphs of girth 5}
	\label{tab:newmoore}
\end{table}

\section{Cayley hypergraphs}\label{sec:cayley}
Our observations above on duality mean that since the dual of a $2$-regular, $d$-uniform hypergraph is simply a graph of the same girth, the cage problems for graphs and hypergraphs have a very natural link. To explore this further, we note that many of the constructions of small graphs of large girth in the survey~\cite{exoo} depend on Cayley graphs or similar constructions. We therefore seek a natural analogue of these constructions in the hypergraph context.

Recall that given a group $G$ and an inverse-closed subset $S$ of $G$, the \emph{Cayley graph} $\Cay(G,S)$ has vertex set the elements of $G$, with an edge from $u$ to $v$ if and only if $u^{-1}v\in S$. There are a number of ways to generalise this idea to the hypergraph setting, but for our purposes the most useful definition is one described in~\cite{Buratti}.
\begin{definition}
Let $G$ be a finite group, $S\subseteq G\setminus\{1\}$ and let $t\geq 2$. The $t$-\emph{Cayley hypergraph} $t$-$\Cay(G,S)$ has vertex set $G$ and hyperedge set

$\quad\{\{g,gs,\ldots,gs^{t-1}\}:g\in G,s\in S\}$.
\end{definition}
Note that taking $t=2$ in the above definition reduces to the usual definition of a Cayley graph.

Our principal interest is in $d$-regular, $r$-uniform hypergraphs, and so we restrict the set $S$ to consist of $d$ elements of order exactly $r$, with suitable conditions on the choice of these elements to ensure that the hypergraph is linear (no two hyperedges sharing more than one vertex). With this restriction, the hyperedges of $r$-$\Cay(G,S)$ are then the left cosets of $\langle s\rangle$ for all $s\in S$.

In the following section we apply this idea to the cage problem for cubic graphs.

\section{Smallest known cubic graphs of given girth}\label{sec:cages}
The problem of finding the smallest 3-regular (cubic) graph of given girth has a long history. For girths of 13 or larger, the optimal value is not known, and for larger girths the smallest known graphs are significantly larger than the Moore bound. (For the history of this problem and the current state of knowledge, see the survey paper~\cite{exoo}.)

It is thus of interest to try to improve the bounds on the smallest cubic graphs of given girth. Much of the focus of previous authors has been on Cayley graphs or other vertex-transitive graphs. It turns out that by searching for ``small'' 2-regular, 3-uniform Cayley hypergraphs of given girth, as described in Section~\ref{sec:cayley}, we can make some improvements to the current record graphs. We note that the graphs constructed by this method are edge-transitive but not necessarily vertex-transitive; thus this population has been less extensively investigated by previous authors.

The method used to find suitable 2-regular, 3-uniform hypergraphs is as follows. The basic idea is to identify a candidate group $G$ generated by two elements $a,b$ both of order 3. We begin by finding all such groups of order no more than $2000$, making use of the library of small groups in \texttt{GAP}~\cite{GAP}. To this list we append all the perfect groups of order up to a million which can be generated by two elements of order 3, again using \texttt{GAP}. (This includes all the simple groups $\PSL(2,q)$ in this range, which have been used by previous authors to good effect.) In addition, we added a number of groups generated by two random elements of order $3$ in a suitable symmetric group. 

We then construct all possible direct products of groups in this list, with the restriction that the resultant group should have order under 2 million and still be generated by two elements of order 3. THis gives a list of $34,970$ candidate groups; of course, this is not an exhaustive list of all possible groups in this range. The upper limit of the orders of groups considered was chosen on grounds of practicality; in fact it was sufficient to find many interesting examples.

For each group $G$ in our candidate list, the algorithm carries out the following steps.

\begin{enumerate}
	\item Find orbit representatives of pairs $a,b$ of elements of order 3 generating $G$. (Or a random sample if there are too many.)
	\item Compute the girth of $H=3$-$\Cay(G,\{a,b\})$. This is the smallest $g$ such that there exists a word $\alpha_1\beta_2\alpha_3\beta_4\cdots\alpha_{g-1}\beta_g=1$, where each $\alpha_i\in\{a,a^{-1}\}$ and $\beta_j\in\{b,b^{-1}\}$.
	\item The dual $H^*$ is then a cubic bipartite graph of order $\frac{2}{3}|G|$ and also has girth $g$.
\end{enumerate}

A table of record smallest cubic graphs for girths up to 32 is maintained on the \texttt{CombinatoricsWiki} website~\cite{wiki}. Our search has resulted in updates to the table of cubic graphs at girths 23, 24, 28, 29, 20, 31 and 32. Table~\ref{tab:cubic} shows the revised table of record graphs. The groups and generating sets giving rise to the 2-regular, 3-uniform hypergraphs used in these constructions are listed in Table~\ref{tab:d2r3} in the Appendix. 

Note that for the odd girths in the table, the graphs have been constructed by excision from a slightly larger graph of girth one greater. Variations on the method of excision have been used by previous authors~\cite{exoo}. Here we have used a simple method based on the ideas of Biggs~\cite{Biggs1998}, where we excise one large tree from the graph using the method of Biggs, then repeatedly excise further 4-vertex trees until no further progress can be made.

\begin{table}[h]\centering
	\begin{tabular}{|crl|}
		\hline
		Girth & Order & Description \\
		\hline
		3 & $4$ & $K_4$ \\
		4 & $6$ & $K_{3,3}$ \\
		5 & $10$ & Petersen \\
		6 & $14$ & Heawood \\
		7 & $24$ & McGee \\
		8 & $30$ & Tutte \\
		9 & $58$ & Brinkmann-McKay-Saager \\
		10 & $70$ & O’Keefe-Wong \\
		11 & $112$ & McKay-Myrvold; Balaban \\
		12 & $126$ & Benson \\
		13 & $272$ & McKay-Myrvold; Hoare \\
		14 & $384$ & McKay; Exoo \\
		15 & $620$ & Biggs \\
		16 & $960$ & Exoo \\
		17 & $2,176$ & Exoo \\
		18 & $2,560$ & Exoo \\
		19 & $4,324$ & Hoare \\
		20 & $5,376$ & Exoo \\
		21 & $16,028$ & Exoo \\
		22 & $16,206$ & Biggs-Hoare \\
		23 & $35,446$ & This paper \\
		24 & $35,640$ & This paper \\
		25 & $108,906$ & Exoo \\
		26 & $109,200$ & Bray-Parker-Rowley \\
		27 & $285,852$ & Bray-Parker-Rowley \\
		28 & $368,640$ & This paper \\
		29 & $805,746$ & This paper \\
		30 & $806,736$ & This paper \\
		31 & $1,440,338$ & This paper \\
		32 & $1,441,440$ & This paper \\
		\hline
	\end{tabular}
	\caption{Smallest known cubic graphs of given girth}
	\label{tab:cubic}
\end{table}

\section{Smallest known graphs and hypergraphs of given girth}\label{sec:hypercages}
As noted in Section~\ref{sec:intro}, a $d$-regular, $r$-uniform hypergraph of girth $g$ can be viewed as an incidence structure having a Levi graph which is a bipartite graph of girth $2g$, and with each vertex in one partition having degree $d$ and each vertex in the other partition degree $r$. Such a graph is called a bipartite $(d,r;2g)$ graph and the study of the smallest such graphs with given parameters was initiated in~\cite{Filipovski}. The natural ``Moore'' bound for such graphs is derived in~\cite{Filipovski} and is essentially identical to the expressions in Equations~\ref{eq:mooreodd1} and~\ref{eq:mooreeven1}.

Here we confine ourselves to tabulating those hypergraphs where the smallest known examples can be determined from the results of our search described in Section~\ref{sec:cages}. A 3-regular, 3-uniform hypergraph of order $n$ and girth $g$ has an incidence graph which is a cubic bipartite graph of order $2n$ and girth $2g$. So the smallest hypergraphs can be determined from the list of smallest known bipartite cubic graphs of even girth. These are tabulated in Table~\ref{tab:d3r3}.

\begin{table}[h]\centering\small
	\setlength{\tabcolsep}{1pt}
	\begin{tabular}{|p{1.5cm}p{2cm}p{2.5cm}|}
		\hline
		Girth & Graph & Hypergraph\\
		\hline
		$6$ & $14$ & $7$ \\
		$8$ & $30$ & $15$ \\
		$10$ & $70$ & $35$ \\
		$12$ & $126$ & $63$ \\
		$14$ & $384$ & $192$ \\
		$16$ & $960$ & $480$ \\
		$18$ & $2,560$ & $1,280$ \\
		\hline
	\end{tabular}
	\quad
	\begin{tabular}{|p{1.5cm}p{2cm}p{2.5cm}|}
		\hline
		Girth & Graph & Hypergraph\\
		\hline
		$20$ & $5,376$ & $2,688$ \\
		$22$ & $16,206$ & $8,103$ \\
		$24$ & $35,640$ & $17,820$ \\
		$26$ & $109,200$ & $54,600$ \\
		$28$ & $368,640$ & $184,320$ \\
		$30$ & $806,736$ & $403,368$ \\
		$32$ & $1,441,440$ & $720,720$ \\
		\hline
	\end{tabular}
	\caption{Smallest known 3-regular, 3-uniform hypergraphs}
	\label{tab:d3r3}
\end{table}

In a similar way, a $2$-regular, $3$-uniform hypergraph can be viewed as the dual of a cubic graph. If this cubic graph had order $n$, then its dual hypergraph has order $3n/2$ and the same girth. This allows us to tabulate these hypergraphs in Table~\ref{tab:d2r3}.

\begin{table}[h]\centering\small
	\setlength{\tabcolsep}{1pt}
	\begin{tabular}{|p{1.5cm}p{2cm}p{2.5cm}|}
		\hline
		Girth & Graph & Hypergraph\\
		\hline
		$3$ & $4$ & $6$ \\
		$4$ & $6$ & $9$ \\
		$5$ & $10$ & $15$ \\
		$6$ & $14$ & $21$ \\
		$7$ & $24$ & $36$ \\
		$8$ & $30$ & $45$ \\
		$9$ & $58$ & $87$ \\
		$10$ & $70$ & $105$ \\
		$11$ & $112$ & $168$ \\
		$12$ & $126$ & $189$ \\
		$13$ & $272$ & $408$ \\
		$14$ & $384$ & $576$ \\
		$15$ & $620$ & $930$ \\
		$16$ & $960$ & $1,440$ \\
		$17$ & $2,176$ & $3,264$ \\
		\hline
	\end{tabular}
	\quad
	\begin{tabular}{|p{1.5cm}p{2cm}p{2.5cm}|}
		\hline
		Girth & Graph & Hypergraph\\
		\hline
		$18$ & $2,560$ & $3,840$ \\
		$19$ & $4,324$ & $6,486$ \\
		$20$ & $5,376$ & $8,064$ \\
		$21$ & $16,028$ & $24,042$ \\
		$22$ & $16,206$ & $24,309$ \\
		$23$ & $35,446$ & $53,169$ \\
		$24$ & $35,640$ & $53,460$ \\
		$25$ & $108,906$ & $163,359$ \\
		$26$ & $109,200$ & $163,800$ \\
		$27$ & $285,852$ & $428,778$ \\
		$28$ & $368,640$ & $552,960$ \\
		$29$ & $805,746$ & $1,208,619$ \\
		$30$ & $806,736$ & $1,210,104$ \\
		$31$ & $1,440,338$ & $2,160,507$ \\
		$32$ & $1,441,440$ & $2,162,160$ \\
		\hline
	\end{tabular}
	\caption{Smallest known 2-regular, 3-uniform hypergraphs}
	\label{tab:d2r3}
\end{table}

\section*{Acknowledgements}
The second author gratefully acknowledges funding support from the London Mathematical Society under reference ECF-2021-27.

\newpage
\section*{Appendix}
Table~\ref{tab:d2r3gen} lists the order, isomorphism class and generators for the groups used to construct the record graphs identified in Table~\ref{tab:cubic}.
\begin{table}[h]\centering\small
    \setlength{\tabcolsep}{0.25em}
	\begin{tabular}{|lll|}
		\hline
		Girth & Order & Group and generators\\
		\hline
		$24$ & $53,460$ & $\PSL(2,11) \times (\Z_3^3 \rtimes \Z_3)$ \\
		& & $(6,9,8)(12,15,17)(13,14,20)(16,18,19)$ \\
		& & $(1,3,6)(2,5,8)(4,7,9)(10,12,16)(11,20,18)(14,15,17)$ \\
		\hline
		$28$ & $552,960$ & $(\Z_2^4 \rtimes \SL(2,5)) \times \SL(2,3) \times \A_4$ \\
		& & $(2,3,4)(6,12,9)(7,8,11)(14,19,18)(15,24,23)(16,17,26)(20,25,21)(22,28,27)$ \\
		& & $\quad (29,30,32)(31,34,33)(35,38,37)(36,40,39)(41,45,44)(42,46,43)(47,50,49)$ \\
		& & $\quad (48,52,51)$\\
		& & $(1,2,3)(5,6,7)(8,9,10)(13,14,15)(16,17,23)(18,27,28)(19,24,21)(22,26,25)$ \\
		& & $\quad (29,34,37)(30,45,43)(31,46,44)(32,39,33)(35,40,49)(36,38,51)(41,52,47)$ \\
		& & $\quad (42,50,48)$ \\
		\hline
		$30$ & $1,210,104$ & $(\Z_7^3 \cdot \PSL(3,2)) \times (\Z_7 \rtimes \Z_3)$ \\
		& & $(2,5,6)(3,4,7)(10,27,28)(11,32,16)(12,24,21)(13,33,14)(15,40,49)(17,23,22)$ \\
		& & $\quad (18,44,38)(19,41,47)(20,35,36)(25,60,54)(26,42,45)(29,61,57)(30,53,52)$ \\
		& & $\quad (31,39,48)(34,43,46)(37,62,55)(50,59,58)(51,63,56)$ \\
		& & $(1,2,3)(4,5,6)(8,11,61)(9,44,33)(10,19,34)(12,17,63)(13,23,53)(14,16,62)$ \\
		& & $\quad (15,26,31)(18,52,40)(21,22,60)(24,32,59)(25,55,41)(27,29,38)(28,45,36)$ \\
		& & $\quad (30,54,43)(35,37,56)(39,50,57)(42,51,58)(46,48,49)$ \\
		\hline
		$32$ & $2,162,160$ & $\Z_3 \times \PSL(2,13) \times \PSL(2,11)$ \\
		& & $(1,2,3)(4,6,10)(5,14,12)(8,9,11)(16,20,17)(18,26,22)(19,23,27)(21,24,28)$ \\
		& & $(6,9,11)(7,8,14)(10,12,13)(15,17,27)(16,26,20)(19,24,21)(22,28,25)$ \\
		\hline
	\end{tabular}
	\caption{Generators for record 2-regular, 3-uniform hypergraphs}
	\label{tab:d2r3gen}
\end{table}

\end{document}